\documentclass[12pt]{amsart}

\theoremstyle{plain}
\newtheorem{theorem}                 {Theorem}      [section]
\newtheorem{proposition}  [theorem]  {Proposition}

\newtheorem{lemma}        [theorem]  {Lemma}

\theoremstyle{definition}

\title[The holomorphic sectional Curvature of General Natural K\" ahler...]{THE HOLOMORPHIC SECTIONAL CURVATURE OF GENERAL NATURAL K\"AHLER STRUCTURES ON COTANGENT BUNDLES}

\author{S.~L.~Dru\c t\u a*}

\thanks{* Partially supported by the Grant TD-158/2007, CNCSIS,
Ministerul  Educa\c tiei \c si Cercet\u arii, Rom\^ania}

\begin{document}

\Large

\maketitle{\footnotesize{\bf Abstract.} We study the conditions
under which a K\"ahlerian structure $(G,J)$ of general natural
lift type on the cotangent bundle $T^*M$ of a Riemannian manifold
$(M,g)$ has constant holomorphic sectional curvature. We obtain
that a certain parameter involved in the condition for
$(T^*M,G,J)$ to be a K\"ahlerian manifold, is expressed as a
rational function of the other two, their derivatives, the
constant sectional curvature of the base manifold $(M,g)$, and the
constant holomorphic sectional curvature of the general natural
K\"ahlerian structure $(G,J)$.


\par}

\normalsize \section{Introduction}

The natural lifts introduced on the cotangent bundle of a
Riemannian ma-\\nifold $(M,g)$ (see \cite{Kolar},
\cite{KowalskiSek}), leaded to some geometric structures studied
in the last years in papers like \cite{OprPap1}, \cite{OprPap2},
\cite{OprPap3}, \cite{OprPor1}, \cite{Porosniuc5},
\cite{Porosniuc3}, \cite{Porosniuc1}, \cite{Porosniuc2}. A part of
the results obtained in these works are similar to some results
from the geometry of the tangent bundle $TM$, the dual of the
cotangent bundle $T^*M$. The differences which appear are related
to the construction of the lifts on the cotangent bundle, the
technics being different from those used in the geometry of the
tangent bundle (see \cite{YanoIsh}).

In the paper \cite{Oproiu4}, Oproiu introduced the general
expression for the natural $1$-st order almost complex structure
$J$ on the tangent bundle $TM$ and the notion of general natural
lifted metric $G$ on $TM$, defined by the Riemannian metric $g$,
with respect to which the horizontal and vertical distributions
are no more orthogonal to each other, contrary to the diagonal
case, treated in \cite{Oproiu3}. The author obtained that the
family of K\" ahlerian structures $(G,J)$ of general natural lift
type on $TM$ depends on three essential parameters (two of them
are involved in the expression of the integrable almost complex
structure $J$, and the third one is a certain proportionality
factor, from the condition for $(G,J)$ to be almost Hermitian).

The present author defined in the paper \cite{Druta}, an almost
complex structure of general natural lifted type on the cotangent
bundle $T^*M$, and a general natural lifted metric to $T^*M$,
obtained from the Riemannian metric $g$ of the base manifold $M$.
The main result is that the family of general natural K\" ahler
structures on $T^*M$ depends on three essential parameters (one is
a certain proportionality factor obtained from the condition for
the structure to be almost Hermitian and the other two are
coefficients involved in the definition of the integrable almost
complex structure $J$ on $T^*M$).

In the joint work \cite{OprDruta}, Oproiu and the present author
studied the conditions under which the K\" ahlerian manifold
$(TM,G,J)$ of general natural lift type has constant holomorphic
sectional curvature. They obtained that the proportionality factor
involved in the condition for $(TM,G,J)$ to be K\" ahlerian is
expressed as a rational function of the two essential parameters
involved in the expression of $J$ (integrable almost complex
structure on $TM$), their derivatives, the constant sectional
curvature of $(M,g)$ and the constant holomorphic sectional
curvature of $(TM,G,J)$.

In the present paper we are interested in finding some properties
of the curvature tensor field $K$ of the general natural K\" ahler
structure $(G,J)$ on the cotangent bundle $T^*M$. Namely, we find
the conditions under which the K\" ahlerian structure considered
on $T^*M$ has constant holomorphic sectional curvature. By doing
some quite long computations with the RICCI package from
Mathematica, we get the expressions of the components of the
curvature tensor field of the manifold $(T^*M,G)$ and those of the
curvature tensor field $K_0$ of the K\" ahlerian manifold
$(T^*M,G,J)$ having constant holomorphic sectional curvature $k$.
The vanishing conditions for the components of difference $K-K_0$
lead to the conclusion that $(T^*M,G,J)$ has constant holomorphic
sectional curvature $k$, if and only if the proportionality factor
involved in the condition for $(T^*M,G,J)$ to be K\" ahlerian is a
rational function depending on the two essential parameters
involved in the expression of the integrable almost complex
structure $J$, their derivatives, the constant sectional curvature
of $(M,g)$ and $k$.

The manifolds, tensor fields and other geometric objects
considered in this paper are assumed to be differentiable of class
$C^\infty $ (i.e. smooth). The Einstein summation convention is
used throughout this paper, the range of the indices
$h,i,j,k,l,m,r $ being always $\{1,\dots ,n\}$.

The author wants to express acknowledgements to Professor Oproiu
for the techniques learned during the elaboration of the joint
work \cite{OprDruta}, for the suggestions, encouragements, and
support throughout the present work and throughout the PhD period.

\vspace*{-.3cm}

\section{Preliminary results}If $(M,g)$ is a smooth
Riemannian manifold of the dimension $n$, and $\pi
:T^*M\rightarrow M$ its cotangent bundle, then the total space
$T^*M$ may be endowed with a structure of a $2n$-dimensional
smooth manifold, induced from the structure of the base manifold,
as follows: from every local chart $(U,\varphi )=(U,x^1,\dots
,x^n)$, it is induced a local chart, $(\pi^{-1}(U),\Phi
)=(\pi^{-1}(U),q^1,\dots , q^n,$ $p_1,\dots ,p_n)$  on $T^*M$,
such that for a cotangent vector $p\in \pi^{-1}(U)\subset T^*M$,
the first $n$ local coordinates $q^1,\dots ,q^n$ are  the local
coordinates of its base point $x=\pi (p)$ in the local chart
$(U,\varphi )$ (in fact we have $q^i=\pi ^* x^i=x^i\circ \pi, \
i=1,\dots n)$; the last $n$ local coordinates $p_1,\dots ,p_n$ of
$p\in \pi^{-1}(U)$ are the vector space coordinates of $p$ with
respect to the natural basis $(dx^1_{\pi(p)},\dots ,
dx^n_{\pi(p)})$, defined by the local chart $(U,\varphi )$,\ i.e.
$p=p_idx^i_{\pi(p)}$.

The notion of $M$-tensor field on the tangent bundle was
introduced in the paper \cite{Mok}. On the cotangent bundle
$T^*M$, an $M$-tensor field of type $(r,s)$ is defined by sets of
$n^{r+s}$ components (functions depending on $q^i$ and $p_i$),
with $r$ upper indices and $s$ lower indices, assigned to induced
local charts $(\pi^{-1}(U),\Phi )$ on $T^*M$, such that the local
coordinate change rule is that of the local coordinate components
of a tensor field of type $(r,s)$ on the base manifold $M$. An
usual tensor field of type $(r,s)$ on $M$ may be thought as an
$M$-tensor field of type $(r,s)$ on $T^*M$. If the considered
tensor field on $M$ is covariant only, the corresponding
$M$-tensor field on $T^*M$ may be identified with the induced
(pullback by $\pi $) tensor field on $T^*M$.

Some useful $M$-tensor fields on $T^*M$ may be obtained as
follows. Let $v,w:[0,\infty ) \rightarrow {\bf R}$ be smooth
functions and let $\|p\|^2=g^{-1}_{\pi(p)}(p,p)$ be the square of
the norm of the cotangent vector $p\in \pi^{-1}(U)$ ($g^{-1}$ is
the tensor field of type (2,0) having the components $(g^{kl}(x))$
which are the entries of the inverse of the matrix $(g_{ij}(x))$
defined by the components of $g$ in the local chart $(U,\varphi
)$). The components $vg_{ij}(\pi(p))$, $p_i$, $w(\|p\|^2)p_ip_j $
define respective $M$-tensor fields of types $(0,2)$, $(0,1)$,
$(0,2)$ on $T^*M$. Similarly, the components $vg^{kl}(\pi(p))$,
$g^{0i}=p_hg^{hi}$, $w(\|p\|^2)g^{0k}g^{0l}$ define respective
$M$-tensor fields of type $(2,0)$, $(1,0)$, $(2,0)$ on $T^*M$. Of
course, all the components considered above are in the induced
local chart $(\pi^{-1}(U),\Phi)$.

\vskip3mm

We recall the splitting of the tangent bundle to $T^*M$ into the
vertical distribution $VT^*M= {\rm Ker}\ \pi _*$ and the
horizontal one determined by the Levi Civita connection $\dot
\nabla $ of $g$:
\begin{equation}
TT^*M=VT^*M\oplus HT^*M.
\end{equation}

If $(\pi^{-1}(U),\Phi)=(\pi^{-1}(U),q^1,\dots ,q^n,p_1,\dots ,p_n)$
is a local chart on $T^*M$, induced from the local chart $(U,\varphi
)= (U,x^1,\dots ,x^n)$, the local vector fields
$\frac{\partial}{\partial p_1}, \dots , \frac{\partial}{\partial
p_n}$ on $\pi^{-1}(U)$ define a local frame for $VT^*M$ over $\pi
^{-1}(U)$ and the local vector fields $\frac{\delta}{\delta
q^1},\dots ,\frac{\delta}{\delta q^n}$ define a local frame for
$HT^*M$ over $\pi^{-1}(U)$, where
$$
\frac{\delta}{\delta q^i}=\frac{\partial}{\partial
q^i}+\Gamma^0_{ih} \frac{\partial}{\partial p_h},\ \ \ \Gamma
^0_{ih}=p_k\Gamma ^k_{ih}
 $$
and $\Gamma ^k_{ih}(\pi(p))$ are the Christoffel symbols of $g$.

The set of vector fields $\{\frac{\partial}{\partial p_1},\dots
,\frac{\partial}{\partial p_n}, \frac{\delta}{\delta q^1},\dots
,\frac{\delta}{\delta q^n}\}$ defines a local frame on $T^*M$,
adapted to the direct sum decomposition (1).

We consider
\begin{equation}
t=\frac{1}{2}\|p\|^2=\frac{1}{2}g^{-1}_{\pi(p)}(p,p)=\frac{1}{2}g^{ik}(x)p_ip_k,
\ \ \ p\in \pi^{-1}(U)
\end{equation}
the energy density defined by $g$ in the cotangent vector $p$. We
have $t\in [0,\infty)$ for all $p\in T^*M$.

The computations will be done in local coordinates, using a local
chart $(U,\varphi)$ on $M$ and the induced local chart
$(\pi^{-1}(U),\Phi)$ on $T^*M$.

We shall use the following lemma, which can be proved easily.

\begin{lemma}\label{lema1}
If $n>1$ and $u,v$ are smooth functions on $T^*M$ such that
$$
u g_{ij}+v p_ip_j=0, \quad u g^{ij}+v g^{0i}g^{0j}=0,\quad or
\quad u\delta ^i_j+vg^{0i} p_j=0,
$$
on the domain of any induced local chart on $T^*M$, then $u=0,\
v=0$.
\end{lemma}

\vskip3mm

In the paper \cite{Druta}, the present author considered the real
valued smooth functions $a_1,\ a_2,\ a_3,\ a_4,\ b_1,\ b_2,\ b_3,\
b_4$ on $[0,\infty)\subset {\bf R}$ and studied a general natural
tensor of type $(1,1)$ on $T^*M$, defined by the relations
\begin{equation}\label{Jinvar}
\left\{
\begin{array}{l}
JX^H_p=a_1(t)
(g_X)^V_p+b_1(t)p(X)p_p^V+a_4(t)X_p^H+b_4(t)p(X)(p^\sharp)_p^H,
\\ \mbox{ }  \\
J\theta^V_p=a_3(t)\theta^V_p+b_3(t)g^{-1}_{\pi(p)}
(p,\theta)p_p^V-a_2(t)(\theta^\sharp)_p^H- b_2(t)g^{-1}_{\pi(p)}
(p,\theta)(p^\sharp)_p^H,
\end{array}
\right.
\end{equation}
in every point $p$ of the induced local card $(\pi^{-1}(U),\Phi)$
on $T^*M$, $\forall ~X \in \mathcal{X}(M), \forall~ \theta \in
\Lambda^1 (M)$, where $g_X$ is the 1-form on $M$ defined by
$g_X(Y)=g(X,Y),\ \forall Y\in \mathcal{X}(M)$,
$\theta^\sharp=g^{-1}_\theta$ is a vector field on $M$ defined by
$g(\theta^\sharp,Y)=\theta (Y)~\forall~ Y \in \mathcal{X}(M)$, the
vector $p^\sharp$ is tangent to $M$ in $\pi (p)$, $p^V$ is the
Liouville vector field on $T^*M$ , and $(p^\sharp)^H$ is the
similar horizontal vector field on $T^*M$.

With respect to the adapted frame $\{\frac{\partial}{\partial
p_i},\frac{\delta}{\delta q^j}\}_{i,j=1,\dots,n} $ on $T^*M$, the
expression (\ref{Jinvar}) becomes

\begin{equation}\label{sist3}
\left\{
\begin{array}{l}
J\frac{\delta}{\delta q^i}=a_1(t)g_{ij}\frac{\partial}{\partial
p_j}+ b_1(t)p_iC+a_4(t)\frac{\delta}{\delta
q^i}+b_4(t)p_i\widetilde C,
\\ \mbox{ }  \\
J\frac{\partial}{\partial p_i}=a_3(t)\frac{\partial}{\partial
p_i}+ b_3(t) g^{0i}C-a_2(t)g^{ij}\frac{\delta}{\delta
q^j}-b_2(t)g^{0i}\widetilde C,
\end{array}
\right.
\end{equation}
where $C =p^V$ is the Liouville vector-field on $T^*M$ and
$\widetilde C=(p^\sharp)^H$ is the corresponding horizontal vector
field on $T^*M$.

We can write also
\begin{equation}\label{sist4}
\left\{
\begin{array}{l}
J\frac{\delta}{\delta q^i}=J^{(1)}_{ij}\frac{\partial}{\partial
p_j}+ J4^j_i\frac{\delta}{\delta q^j},
\\ \mbox{ }  \\
J\frac{\partial}{\partial p_i}=J3^i_j\frac{\partial}{\partial
p_j}-J_{(2)}^{ij}\frac{\delta}{\delta q^j},
\end{array}
\right.
\end{equation}
where
$$
J^{(1)}_{ij}=a_1(t)g_{ij}+b_1(t)p_ip_j,\quad
J4^j_i=a_4(t)\delta^j_i+b_4(t)g^{0j}p_i
$$
\vskip1mm
$$
\quad J3^i_j=a_3(t)\delta^i_j+b_3(t)g^{0i}p_j,\quad
J_{(2)}^{ij}=a_2(t)g^{ij}+b_2(t)g^{0i}g^{0j}.
$$

\begin{theorem}\label{th4}\rm{(\cite{Druta})}
A natural tensor field $J$ of type $(1,1)$ on $T^*M$ given by
$(\ref{sist3})$ or $(\ref{sist4})$ defines an almost complex
structure on $T^*M$, if and only if $a_4=-a_3, b_4=-b_3$ and the
coefficients $a_1,\ a_2,\ a_3,\ b_1,\ b_2$ and $b_3$ are related
by
\begin{equation}\label{rel4}
a_1a_2=1+a_3^2\ ,\ \ \ (a_1+2tb_1)(a_2+2tb_2)=1+(a_3+2tb_3)^2.
\end{equation}
\end{theorem}

\vskip 3mm \bf Remark. \rm From the conditions (\ref{rel4}) we
have that the coefficients $a_1,\ a_2,\ a_1+2tb_1,\ a_2+2tb_2$
have the same sign and cannot vanish. We assume that $a_1>0,\
a_2>0,\ a_1+2tb_1>0,\ a_2+2tb_2>0$ for all $t\geq 0$.

\vskip3mm

\bf Remark. \rm The relations (\ref{rel4}) allow us to express two
of the coefficients $a_1,\ a_2,\ a_3,\ b_1$, $b_2,\ b_3$ as
functions of the other four; e.g. we have:
\begin{equation}\label{inlocuire}
a_2=\frac{1+a_3^2}{a_1},\ \ \
b_2=\frac{2a_3b_3-a_2b_1+2tb_3^2}{a_1+2tb_1}.
\end{equation}

\vskip3mm

The integrability condition for the above almost complex structure
$J$ on a manifold $M$ is characterized by the vanishing of its
Nijenhuis tensor field $N_J$, defined by
$$
N_J(X,Y)=[JX,JY]-J[JX,Y]-J[X,JY]-[X,Y],
$$
for all vector fields $X$ and $Y$ on $M$.

\begin{theorem}\label{th3}\rm{(\cite{Druta})}
Let $(M,g)$ be an $n(>2)$-dimensional connected Riemannian
manifold. The almost complex structure $J$ defined by
{\rm(\ref{sist3})} on $T^*M$ is integrable if and only if $(M,g)$
has constant sectional curvature $c$ and the  coefficients $b_1,\
b_2,\ b_3$ are given by:
\begin{equation}\label{integrab}
\begin{cases}
 b_1=\frac{2 c^2 t a_2^2+2 c t a_1 a_2^\prime+a_1
a_1^\prime -c+3 c a_3^2}{a_1-2 t a_1^\prime-2 c t a_2-4 c t^2
a_2^\prime},\\
 \\
b_2=\frac{2 t a_3^{\prime 2}- 2 t a_1^\prime a_2^\prime+c a_2^2+2
c t a_2 a_2^\prime+a_1 a_2^\prime}{a_1-2 t a_1^\prime-2 c t a_2-4
c t^2 a_2^\prime},\\
\\
 b_3 =\frac{a_1 a_ 3^\prime+ 2 c a_2 a_3+4 c t a_2^\prime a_3-
 2 c t a_ 2 a_3^\prime}{a_1-2 t a_1^\prime-2 c t a_2-4 c t^2 a_2^\prime}.
\end{cases}
\end{equation}
\end{theorem}

\bigskip
\bf Remark. \rm In the diagonal case, where $a_3=0$ it follows
$b_3=0$ too, and we have:
$$
a_2=\frac{1}{a_1},\ b_1=\frac{a_1a_1^\prime-c}{a_1-2ta_1^\prime},\
b_2=\frac{c-a_1a_1^\prime}{a_1(a_1^2-2ct)}.
$$

In the paper cited bellow, the present author defined a Riemannian
metric $G$ of general natural lift type is defined by the
relations
\begin{equation}\label{Ginvar}
\left\{
\begin{array}{l}
G_p(X^H, Y^H) = c_1(t)g_{\pi(p)}(X,Y) + d_1(t)p(X)p(Y),
\\ \mbox{ } \\
G_p(\theta^V,\omega^V) = c_2(t)g^{-1}_{\pi(p)}(\theta,\omega) +
d_2(t)g^{-1}_{\pi(p)}(p,\theta)g^{-1}_{\pi(p)}(p,\omega),
\\ \mbox{ } \\
G_p(X^H,\theta^V) = G_p(\theta^V,X^H)
=c_3(t)\theta(X)+d_3(t)p(X)g^{-1}_{\pi(p)}(p,\theta),
\end{array}
\right.
\end{equation}
$\forall~ X,Y \in \mathcal{X}(M),$ $\forall~ \theta, \omega \in
\Lambda^1(M), \forall~p \in T^*M$.

\vskip3mm

Using the adapted frame $\{\frac{\partial}{\partial
p_i},\frac{\delta}{\delta q^j}\}_{i,j=1,\dots,n} $ on $T^*M$, we
can write the expression (\ref{Ginvar}) in the next form
\begin{equation}\label{rel11}
\left\{
\begin{array}{l}
G(\frac{\delta}{\delta q^i}, \frac{\delta}{\delta
q^j})=c_1(t)g_{ij}+ d_1(t)p_ip_j=G^{(1)}_{ij},
\\   \mbox{ } \\
G(\frac{\partial}{\partial p_i}, \frac{\partial}{\partial
p_j})=c_2(t)g^{ij}+d_2(t)g^{0i}g^{0j}=G_{(2)}^{ij},
\\   \mbox{ } \\
G(\frac{\partial}{\partial p_i},\frac{\delta}{\delta q^j})=
G(\frac{\delta}{\delta q^i},\frac{\partial}{\partial p_j})=
c_3(t)\delta_i^j+d_3(t)p_ig^{0j}=G3_i^j,
\end{array}
\right.
\end{equation}
where $c_1,\ c_2,\ c_3,\ d_1,\ d_2,\ d_3$ are six smooth functions
of the density energy on $T^*M$.

The conditions for $G$ to be positive definite are assured if
\begin{equation}\label{pozdef}
c_1+2td_1>0,\quad c_2+2td_2>0,
\end{equation}
$$
(c_1+2td_1)(c_2+2td_2)-(c_3+2td_3)^2>0.
$$

The metric $G$ is almost Hermitian with respect to the general
almost complex structure $J$, if
$$
G(JX,JY)=G(X,Y),
$$
for all vector fields $X,Y$ on $T^*M$.

The author proved the following result

\begin{theorem}\label{th4}\rm{(\cite{Druta})}{The family of natural, Riemannian metrics  $G$
on $T^*M$ such that $(T^*M,G,J)$ is an almost Hermitian manifold,
is given by {\rm(\ref{rel11})}, provided that the coefficients
$c_1,\ c_2,\  c_3,\  d_1,\ d_2,$ and $d_3$ are related to the
coefficients $a_1,\ a_2,\ a_3,\ b_1,\ b_2,$ and $b_3$ by the
following proportionality relations
\begin{equation}\label{proportiec}
\frac{c_1}{a_1} =\frac{c_2}{a_2}=\frac{c_3}{a_3} = \lambda
\end{equation}
\begin{equation}\label{proportied}
\frac{c_1+2 t d_1}{a_1+2 t b_1} =\frac{c_2+2 t d_2}{a_2+2 t
b_2}=\frac{c_3+2 t d_3}{a_3+2 t b_3} = \lambda +2 t \mu,
\end{equation}
where the proportionality coefficients $\lambda>0 $ and $\lambda
+2t\mu>0$ are functions depending on $t$.}
\end{theorem}

\bf Remark. \rm In the case where $a_3=0$, it follows that
$c_3=d_3=0$ and we obtain the almost Hermitian structure
considered in \cite{OprPor1}, \cite{Porosniuc3}. Moreover, if
$\lambda =1$ and $\mu =0$, we obtain the almost K\"ahlerian
structure considered in the mentioned papers.

\vskip 3mm Considering the two-form $\Omega $ defined by the
almost Hermitian structure $(G,J)$ on $T^*M$
$$
\Omega (X,Y)=G(X,JY),
$$
for all vector fields $X,Y$ on $T^*M$, we obtain the following
result from \cite{Druta}:

\begin{proposition}\label{prop5}\rm{(\cite{Druta})}
The expression of the $2$-form $\Omega $ in the local adapted
frame\\ $\{\frac{\partial}{\partial p_i},\frac{\delta}{\delta
q^j}\}_{i,j=1,\dots,n} $ on $T^*M$, is given by
$$
\Omega \left(\frac{\partial}{\partial
p_i},\frac{\partial}{\partial p_j}\right)=0,\ \Omega \left(
\frac{\delta}{\delta q^i},\frac{\delta}{\delta q^j}\right)=0,\
\Omega\left(\frac{\partial}{\partial p_i},\frac{\delta}{\delta
q^j}\right)= \lambda \delta^i_j+\mu g^{0i}p_j
$$
or, equivalently
\begin{equation}
\Omega =(\lambda \delta^i_j+\mu g^{0i}p_j)Dp_i\wedge dq^j,
\end{equation}
where  $Dp_i=dp_i-\Gamma^0_{ih}dq^h$ is the absolute differential
of $p_i$.
\end{proposition}

Next, by calculating the exterior differential of $\Omega$,  we
may state:
\begin{theorem}\label{th6}\rm{(\cite{Druta})}
The almost Hermitian structure $(T^*M,G,J)$ is almost
K\"{a}hlerian if and only if
$$
\mu=\lambda ^\prime .
$$
\end{theorem}

\vskip 2mm \bf Remark. \rm The family of general almost
K\"ahlerian structures on $T^*M$ depends on five essential
coefficients $a_1,\ a_3,\ b_1,\ b_3,\ \lambda$, which must satisfy
the supplementary conditions $a_1>0,\ a_1+2tb_1>0,\ \lambda>0,\
\lambda +2 t \mu>0$.

The main result obtained in \cite{Druta} is the next one:

\begin{theorem}\rm{(\cite{Druta})}
A general natural lift almost Hermitian structure $(G,J)$ on
$T^*M$ is K\"ahlerian if and only if  the almost complex structure
$J$ is integrable (see Theorem $\ref{th3}$) and
$\mu=\lambda^\prime$.
\end{theorem}

\vskip 2mm

\bf Remark. \rm The family of general natural K\"ahlerian
structures on $T^*M$ depends on three  essential coefficients
$a_1,\ a_3,\ \lambda$, which must satisfy the supplementary
conditions $a_1>0,\ a_1+2tb_1>0,\ \lambda>0$,
$\lambda+2t\lambda^\prime>0$, where $b_1$ is is given by
(\ref{integrab}).

Examples of such structures can be found in \cite{OprPor1},
\cite{Porosniuc3}.

\section{General natural K\"ahler structures of constant
holomorphic sectional curvature on  cotangent bundles}

The Levi-Civita connection $\nabla$ of the Riemannian manifold
$(T^*M,G)$ is obtained from the formula
\begin{align*}
2G(\nabla_XY,Z)&=X(G(X,Z))+Y(G(X,Z))-Z(G(X,Y))\\&
+G([X,Y],Z)-G([X,Z],Y)-G([Y,Z],X); \\&\qquad\qquad\qquad\forall
X,Y,Z\in \chi(M)
\end{align*}
and is characterized by the conditions
$$
\nabla G=0,\ T=0,
$$
where $T$ is the torsion tensor of $\nabla.$

In the case of the cotangent bundle $T^*M$ we can obtain the
explicit expression of $\nabla$. The symmetric $2n\times 2n$
matrix
$$\left(
\begin{matrix}
G^{(1)}_{ij} &  G3^j_i \\
G3^i_j &  G_{(2)}^{ij}
\end{matrix}\right)
$$
associated to the metric $G$ in the base $\{\frac{\delta}{\delta
q^i},\frac{\partial}{\partial p_j}\}_{i,j=1,\dots,n}$ has the
inverse
$$
\left(
\begin{matrix}
H_{(1)}^{ij} &  H3_i^j \\
H3^i_j &  H^{(2)}_{ij}
\end{matrix}
\right)
$$
where the entries are the blocks
\begin{equation}\label{matrinv}
H_{(1)}^{kl}=e_1g^{kl}+f_1g^{0k}g^{0l},\quad
H^{(2)}_{kl}=e_2g_{kl}+f_2p_kp_l,\quad
H3^k_l=e_3\delta^k_l+f_3g^{0k}p_l.
\end{equation}
Here $g^{kl}$ are the components of the inverse of the matrix
$(g_{ij})$, $g^{0k}=p_ig^{ik}$, and $e_1,\ f_1,\ e_2,\ f_2,\ e_3$,
$f_3:[0,\infty)\rightarrow \mathbb{R},$ some real smooth
functions. Their expressions are obtained  by solving the system:
\begin{eqnarray}\label{sistinv}
\begin{cases}
G^{(1)}_{ih}H_{(1)}^{hk}+G3_i^hH3^k_h=\delta_i^k\\
G^{(1)}_{ih}H3^h_k+G3_i^hH^{(2)}_{hk}=0\\
G3^i_hH_{(1)}^{hk}+G_{(2)}^{ih}H3_h^k=0\\
G3^i_hH3^h_k+G_{(2)}^{ih}H^{(2)}_{hk}=\delta^i_k,
\end{cases}
\end{eqnarray}
in which we substitute the relations (\ref{rel11}) and
(\ref{matrinv}). By using Lemma 1,  we get $e_1,\ e_2,\ e_3$ as
functions of $c_1,\ c_2,\ c_3$
\begin{eqnarray}\label{inversa1}
e_1=\frac{c_2}{c_1c_2-c_3^2},\ \ e_2=\frac{c_1}{c_1c_2-c_3^2},\ \
e_3=-\frac{c_3}{c_1c_2-c_3^2}
\end{eqnarray}
and $f_1,\ f_2,\ f_3$ as functions of $c_1,\ c_2,\ c_3,$ $d_1,\
d_2,\ d_3,$ $e_1,\ e_2,\ e_3$

\begin{equation}\label{inversa2}
\begin{aligned}
f_1&=-\frac{c_2d_1e_1 - c_3d_3e_1 - c_3d_2e_3 + c_2d_3e_3 +
2d_1d_2e_1t - 2d_3^2e_1t}{c_1c_2 - c_3^2 + 2c_2d_1t + 2c_1d_2t -
4c_3d_3t + 4d_1d_2t^2 - 4d_3^2t^2},\\\\
f_2&=\frac{(c_3 + 2d_3t)[(d_3e_1 + d_2e_3)(c_1 + 2d_1t) - (d_1e_1
+ d_3e_3)(c_3 + 2d_3t)]}{(c_2 + 2d_2t)[(c_1 + 2d_1t)(c_2 + 2d_2t)
- (c_3 + 2d_3t)^2]}
\end{aligned}
\end{equation}

$$
-\frac{d_2e_2 + d_3e_3}{c_2 + 2d_2t} ,
$$
$$
f_3=-\frac{(d_3e_1 + d_2e_3)(c_1 + 2d_1t) - (d_1e_1 + d_3e_3)(c_3
+ 2d_3t)}{(c_1 + 2d_1t)(c_2 + 2d_2t) - (c_3 + 2d_3t)^2}
$$

Next we can obtain the expression of the Levi Civita connection of
the Riemannian metric $G$ on $T^*M$.
\begin{theorem}
The Levi-Civita connection $\nabla$ of\ $G$ has the following
expression in the local adapted frame $\{\frac{\delta}{\delta
q^i},\frac{\partial}{\partial p_j}\}_{i,j=1,\dots,n}$
$$
\begin{cases}
\displaystyle\nabla_{\frac{\partial}{\partial p_i}}
\frac{\partial}{\partial p_j}=Q^{ij}_{\ \
h}\frac{\partial}{\partial
p_h}+\widetilde{Q}^{ijh}\frac{\delta}{\delta q^h},~
\nabla_{\frac{\delta}{\delta q^i}} \frac{\partial}{\partial
p_j}=(-\Gamma^j_{ih}+\widetilde{P}_{i\ \ h}^{\
j})\frac{\partial}{\partial p_h}+P_i^{\ jh}\frac{\delta}{\delta
q^h}\\\displaystyle \nabla_{\frac{\partial}{\partial p_i}}
\frac{\delta}{\delta q^j}=P_j^{\ ih}\frac{\delta}{\delta
q^h}+\widetilde{P}_{j\ h}^{\ i}\frac{\partial}{\partial p_h},~
\nabla_{\frac{\delta}{\delta q^i}} \frac{\delta}{\delta
q^j}=(\Gamma^h_{ij}+\widetilde{S}_{ij}^{\ \
h})\frac{\delta}{\delta p_h}+S_{ijh}\frac{\partial}{\partial p_h},
\end{cases}
$$
where $\Gamma^h_{ij}$ are the Christoffel symbols of the
connection $\dot\nabla$ and $M$-tensor fields appearing as
coefficients in the above expressions are given by
\begin{equation}\label{PQS}
\begin{cases}
Q^{ij}_{\ \
h}=\frac{1}{2}(\partial^iG^{jk}_{(2)}+\partial^jG^{ik}_{(2)}-
\partial^kG^{ij}_{(2)})H^{(2)}_{kh}+\frac{1}{2}(\partial^iG3^j_k+
\partial^j G3^i_k) H3^k_h,\\
\widetilde{Q}^{ijh}=\frac{1}{2}(\partial^iG^{jk}_{(2)}+\partial^jG^{ik}_{(2)}-
\partial^kG^{ij}_{(2)})H3_k^h+\frac{1}{2}(\partial^iG3^j_k+
\partial^jG3^i_k)H_{(1)}^{kh},\\\\
P_j^{\ ih}=\frac{1}{2}(\partial^iG3^k_j-
\partial^kG3^i_j)H3^h_k+\frac{1}{2}(\partial^iG_{jk}^{(1)}-
R^0_{ljk}G_{(2)}^{li})H_{(1)}^{kh},\\
\widetilde{P}_{j\ h}^{\ i}=\frac{1}{2}(\partial^iG3^k_j-
\partial^kG3^i_j)H^{(2)}_{kh}+\frac{1}{2}(\partial^iG_{jk}^{(1)}-
R^0_{ljk}G_{(2)}^{li})H3^k_h,\\\\
S_{ijh}=\frac{1}{2}(R^{0}_{lij}G_{(2)}^{lk}-\partial^kG_{ij}^{(1)})H^{(2)}_{kh}-c_3R^0_{ijk}H3^k_h,\\
\widetilde{S}^{\ \
h}_{ij}=\frac{1}{2}(R^{0}_{lij}G_{(2)}^{lk}-\partial^kG_{ij}^{(1)})H3^h_k-c_3R^0_{ijk}H_{(1)}^{kh},
\end{cases}
\end{equation}
where $R^h_{kij}$ are the components of the curvature tensor field
of the Levi Civita connection $\dot \nabla$ of the base manifold
$(M,g)$.
\end{theorem}

If we replace in (\ref{PQS}) the relations (\ref{rel11}) which
define the metric $G$, the expressios (\ref{matrinv}) for the
inverse matrix $H$ of $G$, and the formulas (\ref{inversa1}),
(\ref{inversa2}) we obtain the detailed expressions of $P_i^{\
jh},\ Q^{ij}_{\ \ h},\ S_{ijh},\ \widetilde{P}_{j\ h}^{\ i},\
\widetilde{Q}^{ijh},\ \widetilde{S}^{\ \ h}_{ij}.$

The curvature tensor field $K$ of the connection $\nabla$ is
defined by
$$
K(X,Y)Z=\nabla_X\nabla_YZ-\nabla_Y\nabla_XZ-\nabla_{[X,Y]}Z,\ \ \
X,Y,Z\in \mathcal{X}(TM).
$$

By using the local adapted frame $\{\frac{\delta}{\delta
q^i},\frac{\partial}{\partial p_j}\}_{i,j=1,\dots
,n}=\{\delta_i,\partial^j\}_{i,j=1,\dots ,n}$ we obtain the
horizontal and vertical components of the curvature tensor field:

$$
K(\delta_i,\delta_j)\delta_k={QQQQ_{ijk}}^h\delta_h+QQQP_{ijkh}\partial^h,
$$
$$
K(\delta_i,\delta_j)\partial^k={QQPQ_{ij}}^{kh}\delta_h+QQPP_{ij\
\ h}^{\ \ k}\partial^h,
$$
$$
K(\partial^i,\partial^j)\delta_k=PPQQ^{ij\ \ h }_{\ \
k}\delta_h+{PPQP^{ij}}_{kh}\partial^h,
$$
$$
K(\partial^i,\partial^j)\partial^k=PPPQ^{ijkh}\delta_h+{PPPP^{ijk}}_h\partial^h,
$$
$$
K(\partial^i,\delta_j)\delta_k=PQQQ^{i\ \ \ h}_{\
jk}\delta_h+{PQQP^i}_{jkh}\partial^h,
$$
$$
K(\partial^i,\delta_j)\partial^k=PQPQ^{i\ \ kh}_{\
j}\delta_h+PQPP^{i\ \ k}_{\ j\ \ h}\partial^h,
$$
\vskip3mm
where the coefficients are the $M$-tensor fields given by
$$
{QQQQ_{ijk}}^h=\widetilde{S}^{\ \ l}_{jk}\widetilde{S}^{\ \
h}_{il}+P^{\ lh}_iS_{jkl} -\widetilde{S}^{\ \
h}_{jl}\widetilde{S}^{\ \ l}_{ik}-P^{\ lh}_jS_{ikl}-R_{lij}^0P^{\
lh}_k+R^h_{kij}
$$
$$
QQQP_{ijkh}=\widetilde{S}^{\ \ l}_{jk}S_{ilh}+\widetilde{P}^ {\
l}_{i\ h}S_{jkl}- \widetilde{S}^{\ \
l}_{ik}S_{jlh}-\widetilde{P}^{\ l}_{j\ h}S_{ikl}-\widetilde{P}^{\
l}_{k\ h}R^0_{lij}
$$

$$
{QQPQ_{ij}}^{kh}=\widetilde{P}^{\ k}_{j\ \ l}P^{\ lh}_i+P_j^{\
kl}\widetilde{S}^{\ \ h}_{il}- \widetilde{P}^{\ k}_{i\ \ l}P_j^{\
lh}-P^{\ kl}_{i}\widetilde{S}^{\ \
h}_{jl}-R^0_{lij}\widetilde{Q}^{lkh},
$$
$$
QQPP_{ij\ \ h}^{\ \ k}\partial^h=\widetilde{P}^{\ k}_{j\ \
l}\widetilde{P}^{\ l}_{i\ \ h}+P_j^{\ kl}S_{ilh}- \widetilde{P}^{\
k}_{i\ \ l}\widetilde{P}^{\ l}_{j\ h}-P^{\
kl}_iS_{jlh}-R^0_{lij}Q_{\ \ h}^{lk}-R^k_{lij},
$$

$$
PPQQ^{ij\ \ h }_{\ \ k}\delta_h=\partial^iP_k^{\
jh}-\partial^jP_k^{\ ih}+\widetilde{P}^ {\ j}_{k\ \
l}\widetilde{Q}^{ilh}+ P_k^{\ jl}P_l^{\ ih}-\widetilde{P}^{\ i}_{
k\ l}\widetilde{Q}^{jlh}-P_k^{\ il}P_l^{\ jh},
$$
$$
{PPQP^{ij}}_{kh}=\partial^i\widetilde{P}^{\ j}_{ k\
h}-\partial^j\widetilde{P}^{\ i}_{k\ h}+\widetilde{P}^{\ j}_{k\
l}Q_{\ \ h}^{il} +P_k^{\ jl}\widetilde{P}^{\ i}_{ l\
h}-\widetilde{P}^{\ i}_{k\ l}Q_{\ \ h}^{jl}-P_k^{\
il}\widetilde{P}^{\ j}_{l\ \ h},
$$

$$
PPPQ^{ijkh}=\partial^i\widetilde{Q}^{jkh}-\partial^j\widetilde{Q}^{ikh}+
Q_{\ \ l}^{jk}\widetilde{Q}^{ilh}+\widetilde{Q}^{jkl}P^{\
ih}_l-Q_{\ \ l}^{ik}\widetilde{Q}^{jlh} -\widetilde{Q}^{ikl}P_l^{\
jh},
$$
$$
{PPPP^{ijk}}_h=\partial^iQ_{\ \ h}^{jk}-\partial^jQ_{\ \
h}^{ik}+Q_{\ \ l}^{jk}Q_{\ \
h}^{il}+\widetilde{Q}^{jkl}\widetilde{P}^{\ i}_{l\ \ h}- Q^{ik}_{\
\ l}Q_{\ \ h}^{jl}-\widetilde{Q}^{ikl}\widetilde{P}^{\ j}_{l\ \
h},
$$

$$
PQQQ^{i\ \ \ h}_{\ jk}\delta_h=\partial^i\widetilde{S}^{\ \
h}_{jk}+S_{jkl}\widetilde{Q}^{ilh}+\widetilde{S}^{\ \
l}_{jk}P_l^{\ ih} -\widetilde{P}^{\ i}_{k\ l}P_j^{\ lh}-P^{\
il}_k\widetilde{S}^{\ \ h}_{jl}
$$
$$
{PQQP^i}_{jkh}=\partial^iS_{jkh}+\widetilde{S}^{\ \ l}_{jk}Q_{\ \
h}^{il}+\widetilde{S}^{\ \ l}_{jk}\widetilde{P}^{\ i}_{l\ h}
-\widetilde{P}_{k\ l}^{\ i}\widetilde{P}^{\ l}_{j\ h}-P_k^{\
il}S_{jlh}
$$

$$
PQPQ^{i\ \ kh}_{\ j}=\partial^iP_j^{\ kh}+\widetilde{P}^{\ k}_{j\
l}\widetilde{Q}^{ilh}+P_j^{\ kl}P_l^{\ ih}- Q_{\ \ l}^{ik}P_j^{\
lh}-\widetilde{Q}^{ikl}\widetilde{S}^{\ \ h}_{jl},
$$
$$
PQPP^{i\ \ k}_{\ j\ \ \ h}=\partial^i\widetilde{P}^{\ k}_{j\
h}+\widetilde{P}^{\ k}_{j\ l}Q_{\ \ h}^{il}+P_j^{\
kl}\widetilde{P}^{\ i}_{l\ h}- Q_{\ \ l}^{ik}\widetilde{P}^{\
l}_{j\ h}-\widetilde{Q}^{ikl}S_{jlh}.
$$

In order to get the final expressions of the above $M$-tensor
fields, we have to compute the first and second order partial
derivatives with respect to the cotangential coordinates $p_i$ of
the usual tensor fields involved in the definition of the
Riemannian metric $G$.

\begin{align*}
\partial^iG^{(1)}_{jk}&=c_1'g^{0i}g_{jk}+d_1'g^{0i}p_jp_k+
d_1 \delta^i_jp_k+d_1 p_j\delta^i_k
\\\partial^iG_{(2)}^{jk}&=c_2'g^{0i}g^{jk}+d_2'g^{0i}g^{0j}g^{0k}+
d_2 g^{ij}g^{0k}+d_2 g^{0j}g^{ik}
\\\partial^iG3^j_k&=c_3'g^{0i}\delta^j_k+d_3'g^{0i}g^{0j}p_k+
d_3 g^{ij}p_k+d_3 g^{0j}\delta^i_k
\end{align*}
\begin{align*}
\partial^i\partial^jG^{(1)}_{kl}&=c_1''g^{0i}g^{0j}g_{kl}+c_1'g^{ij}g_{kl}+
d_1''g^{0i}g^{0j}p_kp_l+ d_1'g^{ij}p_kp_l
\\
&+d_1'g^{0j}\delta^i_kp_l+d_1'g^{0j}p_k\delta^i_l+d_1'g^{0i}\delta^j_kp_l+
d_1 \delta^j_k\delta^i_l
\\
&+d_1'g^{0i}p_k\delta^j_l+d_1 \delta^i_k\delta^j_l
\\\partial^i\partial^jG_{(2)}^{kl}&=c_2''g^{0i}g^{0j}g^{kl}+c_2'g^{ij}g^{kl}+
d_2''g^{0i}g^{0j}g^{0k}g^{0l}+ d_2'g^{ij}g^{0k}g^{0l}
\\
&+d_2'g^{0j}g^{ik}g^{0l}+d_2'g^{0j}g^{0k}g^{il}+d_2'g^{0i}g^{jk}g^{0l}+
d_2 g^{jk}g^{il}
\\
&+d_2'g^{0i}g^{0k}g^{jl}+d_2 g^{ik}g^{jl}
\\
\partial^i\partial^jG3^k_l&=c_3''g^{0i}g^{0j}\delta^k_l+c_3'g^{ij}\delta^k_l+
d_3''g^{0i}g^{0j}g^{0k}p_l+ d_3'g^{ij}g^{0k}p_l
\\
&+d_3'g^{0j}g^{ik}p_l+d_3'g^{0j}g^{0k}\delta^i_l+d_3'g^{0i}g^{jk}p_l+
d_3 g^{jk}\delta^i_l
\\
&+d_3'g^{0i}g^{0k}\delta^j_l+d_3 g^{ik}\delta^j_l
\end{align*}
\begin{align*}
\partial^iH_{(1)}^{jk}&=e_1'g^{0i}g^{jk}+f_1'g^{0i}g^{0j}g^{0k}+
f_1g^{ij}g^{0k}+f_1 g^{0j}g^{ik}
\\
\partial^iH^{(2)}_{jk}&=e_2'g^{0i}g_{jk}+f_2'g^{0i}p_jp_k+
f_2\delta^i_jp_k+f_2 p_j\delta^i_k
\\
\partial^iH3^j_k&=e_3'g^{0i}\delta^j_k+f_3'g^{0i}g^{0j}p_k+
f_3g^{ij}p_k+f_3 g^{0j}\delta^i_k
\end{align*}

We get the first order partial derivatives of the $M$-tensor
fields $P_i^{\ jh},\ Q^{ij}_{\ \ h},$ $ S_{ijh},$
$\widetilde{P}_{j\ h}^{\ i},\ \widetilde{Q}^{ijh},\
\widetilde{S}^{\ \ h}_{ij}$ with respect to the cotangential
coordinates $p_i$  and we replace these derivatives, and the
expressions (\ref{inversa1}), (\ref{inversa2}) of the functions
$e_1,\ e_2,\ e_3,$ $f_1,\ f_2,\ f_3$ and of their derivatives in
order to obtain the components of the curvature tensor as
functions of $a_1,\ a_2,\ a_3$ and their derivatives  of first,
second and  third order only. The expressions are obtained by
using the Mathematica package RICCI.

\begin{align*}
\partial^iQ_{\ \ h}^{jk}&=\frac{1}{2}\partial^iH^{(2)}_{lh}(\partial^jG^{(2)}_{kl}+
\frac{1}{2}H^{(2)}_{lh}(\partial^i\partial^jG^{(2)}_{kl}+\partial^i\partial^kG_{(2)}^{jl}
-\partial^i\partial^lG_{(2)}^{jk})
\\
&+\frac{1}{2}\partial^iH3^l_h(\partial^jG3^k_l+\partial^k G3^k_l)+
\frac{1}{2}H3^l_h(\partial^i\partial^jG3^k_l+\partial^i\partial^k
G3^k_l),
\\
\partial_i\widetilde{Q}^{jkh}&=\frac{1}{2}\partial^iH3_l^h(\partial^jG^{kl}_{(2)}+\partial^kG^{jl}_{(2)}-
\partial^lG^{jk}_{(2)})
\\
&+\frac{1}{2}H3_l^h(\partial^i\partial^jG^{kl}_{(2)}+\partial^i\partial^kG^{jl}_{(2)}
-\partial^i\partial^lG^{jk}_{(2)})
\\
&+\frac{1}{2}\partial^iH_{(1)}^{lh}(\partial^jG3^k_l+
\partial^kG3^j_l)
+\frac{1}{2}H_{(1)}^{lh}(\partial^i\partial^jG3^k_l+\partial^i\partial^kG3^j_l),
\end{align*}
\begin{align*}
\partial^i\widetilde{P}_{j\ h}^{\ k}&=\frac{1}{2}\partial^iH^{(2)}_{lh}(\partial^kG3^l_j-
\partial^lG3^k_j)+ \frac{1}{2}H^{(2)}_{lh}(\partial^i\partial^kG3^l_j
-\partial^i\partial^lG3^k_j)
\\
&+\frac{1}{2}\partial^iH3^l_h(\partial^kG_{jl}^{(1)}-
R^0_{mjl}G_{(2)}^{mk})
\\
&+ \frac{1}{2}H3^l_h(\partial^i\partial^kG_{jl}^{(1)}-
R^i_{mjl}G_{(2)}^{mk}-R^0_{mjl}\partial^iG_{(2)}^{mk}),
\\
\partial^iP_j^{\ kh}&=\frac{1}{2}\partial^iH3^h_l(\partial^kG3^l_j-
\partial^lG3^k_j)+
\frac{1}{2}H3^h_l(\partial^i\partial^kG3^l_j
-\partial^i\partial^lG3^k_j)
\\
&+\frac{1}{2}\partial^iH_{(1)}^{hl}(\partial^kG_{jl}^{(1)}-
R^0_{mjl}G_{(2)}^{mk})
\\
&+\frac{1}{2}H_{(1)}^{hl}(\partial^i\partial^kG_{jl}^{(1)}-
R^i_{mjl}G_{(2)}^{mk}-R^0_{mjl}\partial^iG_{(2)}^{mk}),
\end{align*}
\begin{align*}
\partial^iS_{jkh}&=\frac{1}{2}[(c_2^{'}g^{0i}R^{0}_{mjk}+c_2R^i_{mjk}-\partial^i\partial^lG^{(1)}_{jk})H^{(2)}_{lh}+
(c_2R^{0}_{mjk}-\partial^lG_{jk}^{(1)})\partial^iH^{(2)}_{lh}]
\\
&-c_3^{'}g^{0i}R^0_{jkl}H3^l_h-c_3(R^i_{jkl}H3^l_h+R^0_{jkl}\partial^iH3^l_h),
\\
\partial^i\widetilde{S}^{\ \ h}_{jk}&=\frac{1}{2}[(c_2^{'}g^{0i}R^{0}_{mjk}+c_2R^i_{mjk}-\partial^i\partial^lG^{(1)}_{jk})H3^h_l+
(c_2R^{0}_{mjk}-\partial^lG_{jk}^{(1)})\partial^iH3^h_l]
\\
&-c_3'g^{0i}R^0_{jkl}H_{(1)}^{lh}-c_3(R^i_{jkl}H_{(1)}^{lh}+R^0_{jkl}\partial^iH_{(1)}^{lh}).
\end{align*}

\vskip3mm The  tensor field corresponding to the curvature tensor
field of a K\"alerian manifold $(T^*M,G,J)$ having constant
holomorphic sectional curvature $k$ is given by the formula:
\begin{align*}
K_0(X,Y)Z&=\frac{k}{4}[G(Y,Z)X-G(X,Z)Y + G(JY,Z)JX
\\
&-G(JX,Z)JY+2G(X,JY)JZ]
\end{align*}

With respect to the adapted frame
$\{\delta_i,\partial^j\}_{i,j=1,\dots ,n}$, the expressions are
$$
K_0(\delta_i,\delta_j)\delta_k={{QQQQ_0}_{ijk}}^h\delta_h+{QQQP}_{0ijkh}\partial^h,
$$
$$
K_0(\delta_i,\delta_j)\partial^k={{QQPQ_0}_{ij}}^{kh}\delta_h+{QQPP_0}_{ij\
\ h}^{\ \ k}\partial^h,
$$

$$
K_0(\partial^i,\partial^j )\delta_k={PPQQ_0}^{ij\ \ h}_{\ \
k}\delta_h+{{PPQP_0}^{ij}}_{kh}\partial^h,
$$
$$
K_0(\partial^i,\partial^j)\partial^k={{PPPP_0}^{ijk}}_h\partial^h+{PPPQ_0}^{ijkh}\delta_h,
$$

$$
K_0(\partial^i,\delta_j)\delta_k={PQQQ_0}^{i\ \ \ h}_{\
jk}\delta_h+{{PQQP_0}^i}_{jkh}\partial^h,
$$
$$
\quad K_0(\partial^i,\delta_j)\partial^k={PQPP_0}^{i\ \ k}_{\ j\ \
h}\partial^h+{PQPQ_0}^{i\ \ k \ h}_{\ j}\delta_h,
$$
where
$$
{{QQQQ_0}_{ijk}}^h=\frac{k}{4}[G^{(1)}_{jk}\delta_i^h-G^{(1)}_{ik}\delta^h_j-
J3^{h}_i(J^{(1)}_{jl}G3^l_k-J3^l_jG^{(1)}_{lk})+
$$
$$
J3^h_j(J^{(1)}_{il}G3^l_k-J3^l_iG^{(1)}_{lk})
-2J3^h_k(J^{(1)}_{jl}G3^l_i-J3^l_jG^{(1)}_{il})],
$$
$$
{QQQP_0}_{ijkh}=\frac{k}{4}[J^{(1)}_{ih}(J^{(1)}_{jl}G3^l_k-J3^l_jG^{(1)}_{lk})-J^{(1)}_{jh}(J^{(1)}_{il}G3^l_k-J3^l_iG^{(1)}_{lk})+
$$
$$
2J^{(1)}_{kh}(J^{(1)}_{jl}G3^l_i-J3^l_jG^{(1)}_{il})],
$$
$$
{{QQPQ_0}_{ij}}^{kh}=\frac{k}{4}[G3^k_j\delta_i^h-G3^k_i\delta^h_j-
J3^h_i(J^{(1)}_{jl}G_{(2)}^{lk}-J3^l_jG3^k_l)+
$$
$$
J3^h_j(J^{(1)}_{il}G_{(2)}^{lk}-J3^l_iG3^k_l)
-2J_{(2)}^{kh}(J^{(1)}_{jl}G3^l_i-J3^l_jG^{(1)}_{il})],
$$
$$
{QQPP_0}_{ij\ \ \ h}^{\ \ \
k}=\frac{k}{4}[J^{(1)}_{ih}(J^{(1)}_{jl}G_{(2)}^{lk}-J3^l_jG3^k_l)-
J^{(1)}_{jh}(J^{(1)}_{il}G_{(2)}^{lk}-J3^l_iG3^k_l)+
$$
$$
2J3^k_h(J^{(1)}_{jl}G3^l_i-J3^l_jG^{(1)}_{il})],
$$
$$
{PPQQ_0}^{ij\ \ h}_{\ \
k}=\frac{k}{4}[-J_{(2)}^{ih}(J3^j_lG3^l_k-J_{(2)}^{jl}G^{(1)}_{lk})+
J_{(2)}^{jh}(J3^i_lG3^l_k-J_{(2)}^{il}G^{(1)}_{lk})-
$$
$$
2J3_k^h(J3^j_lG_{(2)}^{il}-J_{(2)}^{jl}G3^i_l)],
$$
$$
{{PPQP_0}^{ij}}_{kh}=\frac{k}{4}[G3^j_k\delta_h^ih-G3^i_k\delta^j_h+
J3^i_h(J3^j_lG3^l_k-J_{(2)}^{jl}G^{(1)}_{lk})-
$$
$$
J3^j_h(J3^i_lG3^l_k-J_{(2)}^{il}G^{(1)}_{lk})+
2J^{(1)}_{kh}(J3^j_lG_{(2)}^{il}-J_{(2)}^{jl}G3^i_l)],
$$
$$
{PPPQ_0}^{ijkh}=\frac{k}{4}[-J_{(2)}^{ih}(J3^j_lG_{(2)}^{lk}-J_{(2)}^{jl}G3^k_l)+
J_{(2)}^{jh}(J3^i_lG_{(2)}^{lk}-J_{(2)}^{il}G3^k_l)-
$$
$$
2J_{(2)}^{kh}(J3^j_lG_{(2)}^{il}-J_{(2)}^{jl}G3^i_l)],
$$
$$
{{PPPP_0}^{ijk}}_h=\frac{k}{4}[G_{(2)}^{jk}\delta_h^i-G_{(2)}^{ik}\delta^j_h+
J3_h^i(J3^j_lG_{(2)}^{lk}-J_{(2)}^{jl}G3^k_l)-
$$
$$
J3_h^j(J3^i_lG_{(2)}^{lk}-J_{(2)}^{il}G3^k_l)+
2J3^k_h(J3^j_lG_{(2)}^{il}-J_{(2)}^{jl}G3^i_l)],
$$
$$
{PQQQ_0}^{i\ \ \ h}_{\ jk}=\frac{k}{4}[-G3^i_k\delta_j^h-
J_{(2)}^{ih}(J^{(1)}_{jl}G3^l_k-J3^l_jG^{(1)}_{lk})+
$$
$$
J3^h_j(J3_l^iG3_k^l-J_{(2)}^{il}G^{(1)}_{lk})-
2J3^h_k(J^{(1)}_{jl}G_{(2)}^{il}-J3^l_jG3^i_l)],
$$
$$
{{PQQP_0}^i}_{jkh}=\frac{k}{4}[G^{(1)}_{jk}\delta_h^i+
J3^i_h(J^{(1)}_{jl}G3^l_k-J3^l_jG^{(1)}_{lk})-
$$
$$
J^{(1)}_{jh}(J3_l^iG3_k^l-J_{(2)}^{il}G^{(1)}_{lk})+
2J^{(1)}_{kh}(J^{(1)}_{jl}G_{(2)}^{il}-J3^l_jG3^i_l)],
$$
$$
{PQPQ_0}^{i\ \ kh}_{\
j}=\frac{k}{4}[-G_{(2)}^{ik}\delta^h_j-J_{(2)}^{ih}(J^{(1)}_{jl}G_{(2)}^{lk}-J3^l_jG3^k_l)+
$$
$$
J3^h_j(J3^i_lG_{(2)}^{lk}-J_{(2)}^{il}G3^k_l)-
2J_{(2)}^{kh}(J^{(1)}_{jl}G_{(2)}^{il}-J3^l_jG3^i_l)],
$$
$$
{PQPP_0}^{i\ \ k}_{\ j\ \
h}=\frac{k}{4}[G3_j^k\delta_h^i+J3^i_h(J^{(1)}_{jl}G_{(2)}^{lk}-J3^l_jG3^k_l)-
$$
$$
J^{(1)}_{jh}(J3^i_lG_{(2)}^{lk}-J_{(2)}^{il}G3^k_l)+
2J3^k_h(J^{(1)}_{jl}G_{(2)}^{il}-J3^l_jG3^i_l)].
$$

\vskip2mm

The K\"ahlerian manifold $(T^*M,G,J)$ is of constant holomorphic
sectional curvature if and only if all the components of the
difference $K-K_0$ vanish. In the study of the vanishing
conditions for the components of $K-K_0$ we use following result
similar to the lemma \ref{lema1}.

\begin{lemma}\label{lema2}
If $\alpha _1,\dots , \alpha_{10}$ are smooth functions on $T^*M$
such that
\begin{equation}\label{expr}
\alpha_1g_{hj}g^{ik}+\alpha_2\delta^i_h\delta^k_j+\alpha_3\delta^k_h\delta^i_j+\alpha_4g^{ik}p_hp_j+
\alpha_5\delta^k_jp_hg^{0i}+\alpha_6\delta^k_hp_jg^{0i}+
\end{equation}
$$
\alpha_7\delta^i_jp_hg^{0k}+ \alpha_8g_{hj}g^{0i}g^{0k}
+\alpha_9\delta^i_hp_jg^{0k}+\alpha_{10}p_hp_jg^{0i}g^{0k}=0,
$$
then $\alpha_1=\dots =\alpha_{10}=0$.
\end{lemma}

\emph{Proof:} If we multiply the expression (\ref{expr}) by
$g^{hj}g_{ik}$, we have
$$
\alpha_1 n^2 +\alpha_2 n + \alpha_3 n + 2\alpha_4 n t + 2 \alpha_5
t + 2 \alpha_6 t + 2 \alpha_7 t + 2 \alpha_8 n t  + 2 \alpha_9 t +
4 \alpha_{10} t^2=0.
$$

Since the expression does not depend on the dimension $n$ of the
base manifold, we obtain that
$$
\alpha_1=0,\ \alpha_2 + \alpha_3  + 2 (\alpha_4 + \alpha_8) t = 0,
\  (\alpha_5 + \alpha_6 + \alpha_7 + \alpha_9)t+ 2 \alpha_{10}
t^2=0.
$$

Similarly, we get that $\alpha _2$ and $\alpha_3$ are also zero,
if we multiply the expression (\ref{expr}), respectively by
$\delta_i^h\delta_k^j$ and $\delta_k^h\delta_i^j$.

The product between (\ref{expr}) and $g_{ik}g^{0h}g^{0j},\
\delta_k^jg^{0h}p_i,\ \delta_k^hg^{0j}p_i, \ \delta_i^jg^{0h}p_k,
g^{hj}p_ip_k,$ or $\delta_i^hg^{0j}p_k,$  leads to some expression
in which the coefficients of $n$ are, respectively $2(\alpha_1 t+
2 \alpha_4 t^2)$, $2(\alpha_2 t+ 2 \alpha_5 t^2)$, $2(\alpha_3 t+
2 \alpha_6 t^2)$, $2(\alpha_1 t+ 2 \alpha_8 t^2)$, $2(\alpha_3 t+
2 \alpha_7 t^2)$, $2(\alpha_1 t+ 2 \alpha_8 t^2)$, $2(\alpha_2 t+
2 \alpha_9 t^2)$. This expressions must vanish for all $t\geq 0$.
Since $\alpha_1=\alpha_2=\alpha_3=0$, we obtain that
$\alpha_4=\dots=\alpha_9=0$ too.

Multiplying by $g^{0h}g^{0j}p_ip_k$, the relation (\ref{expr})
becomes
\begin{equation}\label{alpha10}
4[(\alpha_1+\alpha_2+\alpha_3)t^2+2(\alpha_4+\alpha_5+\alpha_6+\alpha_7+\alpha_8+\alpha_9)t^3+4\alpha_{10}t^4]=0.
\end{equation}

Taking into account that $\alpha_i=0,\ \forall i=1,\dots,9$, it
follows from (\ref{alpha10}) that $\alpha_{10}=0.$

\vskip3mm

The final theorem gives the condition under which the K\"ahlerian
manifold of general natural lift type has constant holomorphic
sectional curvature

\begin{theorem}\label{teoremafinala}
The K\"ahlerian manifold $(T^*M,G,J)$ with $G$ and $J$ obtained as
natural lifts of general type of the Riemannian metric $g$ on the
Riemannian manifold $(M,g)$, has constant holomorphic sectional
curvature $k$ if and only if the parameter $\lambda$ is expressed
by
\begin{eqnarray}\label{lambda}
\lambda=\frac{4a_1c}{k(a_1^2 + 2ct + 2a_3^2ct)}
\end{eqnarray}
\end{theorem}

\vskip3mm \emph{Proof.} The expressions of the differences that we
study are quite long, but in  $PQPP^{i\ \ k}_{\ j\ \
h}-{PQPP_0}^{i\ \ k}_{\ j\ \ h}$ two coefficients have shorter
expressions. From the first term which contains $g_{hj}g^{ik}$, by
imposing the annulation of the coefficient, we get

\begin{eqnarray}\label{lambdapr}
\lambda'=-\lambda\frac{a_1'(a_1^2- 2ct - 2a_3^2ct) + 2a_1c(1 +
a_3^2 + 2a_3a_3't)}{a_1(a_1^2 + 2ct + 2a_3^2ct)},
\end{eqnarray}

If we substitute this expression in the second term (which
contains $\delta^i_h \delta^k_j$) we obtain the value of $\lambda$
given by (\ref{lambda}).

The expression of $\lambda'$ obtained by differentiating the
relation (\ref{lambda}), coincides with that obtained by replacing
$\lambda$ in (\ref{lambdapr}). Using RICCI, we prove that all the
components of the difference $K-K_0$ are zero, when the obtained
values of $\lambda^\prime,\ \lambda^{\prime \prime}$ and
$\lambda^{\prime \prime \prime}$ are replaced in these components.
The computation of some differences, such as $PQPP^{i\ \ k}_{\ j\
\ h}-{PQPP_0}^{i\ \ k}_{\ j\ \ h}$,$PQPQ^{i\ \ kh}_{\
j}-{PQPQ_0}^{i\ \ kh}_{\ j}$,
${{PQQP}^i}_{jkh}-{{PQQP_0}^i}_{jkh}$, and $PQQQ^{i\ \ \ h}_{\
jk}-{PQQQ_0}^{i\ \ \ h}_{\ jk}$ is quite hard, since after
imposing the integrability conditions for the almost complex
structure $J$, the expressions become very long, and the command
TensorSimplify did not work on a PC with a RAM memory of 2GB. Thus
I had to impose the integrability conditions in every coefficient
appearing  in the above differences, and to sum the expressions
afterwards.

\bf{Remark}. \rm If $a_3=0$ we obtain the condition for
$(T^*M,G,J)$ to have constant holomorphic sectional curvature in
the case where $G,J$ are natural lifts of diagonal type (see
\cite{Oproiu4}, \cite{OprPor1}).

\strut\hfill simonadruta@yahoo.com
\begin{flushright}
Faculty of Mathematics, \\
University "Al.I. Cuza", \\
Ia\c si 700 506 ROM\^{A}NIA
\end{flushright}


\begin{thebibliography}{11}

\bibitem{Druta}{\bf Dru\c t\u a, S.L.}, {\it Cotangent Bundles with General Natural K\" ahler
structures}. Accepted to publish in R$\acute{e}$vue Roumaine de
Math$\acute{e}$matiques Pures et Appliqu$\acute{e}$es.

\bibitem{OprDruta}{\bf  Druta, S., Oproiu, V.,} {\it General Natural K\"ahler Structures
of Constant Holomorphic Sectional Curvature on Tangent Bundles},
An.St.Univ. "Al.I.Cuza" Iasi (S. N.) Matematica, Tom LIII, 2007,
f.1, 149-166.

\bibitem{Kolar}{\bf Kol\'a\v r, I.; Michor, P.; Slovak, J.}, {\it Natural
Operations in Differential Geometry}, Springer Verlag, Berlin,
1993, vi, 434 pp.

\bibitem{KowalskiSek} {\bf Kowalski, O.; Sekizawa, M.}, {\it Natural
transformations of Riemannian metrics on manifolds to metrics on
tangent bundles-a classification,} Bull. Tokyo Gakugei Univ. (4),
40 (1988), 1-29.

\bibitem{Mok}{\bf Mok, K.P.; Patterson, E.M.; Wong, Y.C.} -- {\it Structure
of symmetric tensors of type (0,2) and tensors of type (1,1) on
the tangent bundle},~ Trans. Amer. Math. Soc.,~234
(1977),~253-278.

\bibitem{Oproiu4}{\bf Oproiu, V.}, {\it A Generalization of Natural Almost Hermitian
Structures on the Tangent Bundles}, Math. J. Toyama Univ.
\textbf{22} (1999) 1--14.

\bibitem{Oproiu3}{\bf Oproiu, V.} -- {\it Some new geometric structures on the
tangent bundles}, Publ. Math. Debrecen, 55/3-4 (1999), 261-281.

\bibitem{OprPap1} {\bf Oproiu, V.; Papaghiuc, N.}, {\it A pseudo-Riemannian
structure on the cotangent bundle}, An. \c St. Univ. "Al. I. Cuza"
Ia\c si, 36 (1990), 265-276.

\bibitem{OprPap2} {\bf Oproiu, V.; Papaghiuc, N.}, {\it Another pseudo-Riemannian
structure on the cotangent bundle}, Blt. Instit. Politehn. Ia\c
si, Tomul 37 (41), Fasc. 1-4, 1991, Sec\c t. I, 27-33.

\bibitem{OprPap3} {\bf Oproiu, V.; Papaghiuc, N.}, {\it Locally
symmetric cotangent bundles}, Matematicki Vesnik, 42 (1990),
221-232.

\bibitem{OprPor1} {\bf Oproiu, V.; Poro\c sniuc, D.D.}, {\it A K\"ahler
Einstein structure on the cotangent bundle of a Riemannian
manifold}, An. \c Stiin\c t. Univ. Al. I. Cuza, Ia\c si 49, s.I,
Mathematics, 2003, f.2, 399-414.

\bibitem{Porosniuc5}{\bf Poro\c sniuc, D.D.}, {\it A class of locally symmetric
K\"ahler Einstein structures on the nonzero cotangent bundle of a
space form}, Balkan Journal of Geometry and Its Applications, 9
(2004), no.2, 68-81.

\bibitem{Porosniuc3}{\bf Poro\c sniuc, D.D.}, {\it A K\"ahler Einstein structure
 on the nonzero cotangent bundle of a space
form},~Italian Journal of Pure and Applied Mathematics,
Udine,~no.18 (2005), 223-234.

\bibitem{Porosniuc1}{\bf Poro\c sniuc, D.D.}, {\it A locally symmetric
K\"ahler Einstein structure on the cotangent bundle of a space
form},~Balkan Journal of Geometry and Its Applications,~9 (2004),
no.1, 94-103.

\bibitem{Porosniuc2}{\bf Poro\c sniuc, D.D.}, {\it A locally symmetric
K\"ahler Einstein structure on a tube in the nonzero cotangent
bundle of a space form}, An. \c Stiin\c t. Univ. Al. I. Cuza, Ia\c
si, 50 (2004), No. 2, 315-326.

\bibitem{YanoIsh}{\bf Yano, K., Ishihara, S.}, {\it Tangent and Cotangent
Bundles}, M. Dekker Inc., New York, 1973.

\end{thebibliography}
\end{document}